\documentclass[11pt]{amsart}
\usepackage{amssymb}
\usepackage{amsmath}
\usepackage[active]{srcltx}
\usepackage{t1enc}
\usepackage[latin2]{inputenc}
\usepackage{verbatim}
\usepackage{amsmath,amsfonts,amssymb,amsthm}
\usepackage[mathcal]{eucal}
\usepackage{enumerate}
\usepackage[centertags]{amsmath}
\usepackage{graphics}

\newtheorem{theorem}{Theorem}

\newtheorem{lemma}{Lemma}

\begin{document}
\author{George Tephnadze}
\title[Fejér means]{Strong convergence theorems of Walsh-Fejér means }
\address{G. Tephnadze, Department of Mathematics, Faculty of Exact and
Natural Sciences, Tbilisi State University, Chavchavadze str. 1, Tbilisi
0128, Georgia}
\email{giorgitephnadze@gmail.com}
\date{}
\maketitle

\begin{abstract}
As main result we prove that Fejér means of Walsh-Fourier series are
uniformly bounded operators from $\ H_{p}$ to $H_{p}$ $\left( 0<p\leq
1/2\right) .$
\end{abstract}

\date{}



\section{INTRODUCTION}

It is well-known that Walsh-Paley system forms not basis in the space $%
L_{1}\left( G\right) .$ Moreover, there is a function in the dyadic Hardy
space $H_{1}\left( G\right) ,$ such that the partial sums of $F$ \ are not
bounded in $L_{1}$-norm. However, in Simon \cite{Si3} the following
estimation was obtained for all $F\in H_{1}\left( G\right) :$%
\begin{equation*}
\frac{1}{\log n}\overset{n}{\underset{k=1}{\sum }}\frac{\left\Vert
S_{k}F\right\Vert _{1}}{k}\leq c\left\Vert F\right\Vert _{H_{1}},
\end{equation*}%
where $S_{k}F$ denotes the $k$-th partial sum of the Walsh-Fourier series of
$F.$ (For the trigonometric analogue see in Smith \cite{sm}, for the
Vilenkin system in Gát \cite{gat1}). Simon \cite{si1} (see also \cite{sw}
and \cite{We}) proved that there is an absolute constant $c_{p},$ depending
only on $p,$ such that
\begin{equation}
\frac{1}{\log ^{\left[ p\right] }n}\overset{n}{\underset{k=1}{\sum }}\frac{%
\left\Vert S_{k}F\right\Vert _{p}^{p}}{k^{2-p}}\leq c_{p}\left\Vert
F\right\Vert _{H_{p}}^{p},\text{ \ \ }\left( 0<p\leq 1\right) ,  \label{1cc}
\end{equation}%
for all $F\in H_{p}$, where $\left[ p\right] $ denotes integer part of $p.$

Weisz \cite{We3} considered the norm convergence of Fejér means of
Walsh-Fourier series and proved the following:

\textbf{Theorem W1 (Weisz).} Let $p>1/2$ and $F\in H_{p}.$ Then
\begin{equation*}
\left\Vert \sigma _{k}F\right\Vert _{H_{p}}\leq c_{p}\left\Vert F\right\Vert
_{H_{p}}.
\end{equation*}

This theorem implies that%
\begin{equation*}
\frac{1}{n^{2p-1}}\overset{n}{\underset{k=1}{\sum }}\frac{\left\Vert \sigma
_{k}F\right\Vert _{H_{p}}^{p}}{k^{2-2p}}\leq c_{p}\left\Vert F\right\Vert
_{H_{p}}^{p},\text{ \ \ \ }\left( 1/2<p<\infty \right) .
\end{equation*}

If Theorem W1 holds for $0<p\leq 1/2,$ then we would have%
\begin{equation}
\frac{1}{\log ^{\left[ 1/2+p\right] }n}\overset{n}{\underset{k=1}{\sum }}%
\frac{\left\Vert \sigma _{k}F\right\Vert _{H_{p}}^{p}}{k^{2-2p}}\leq
c_{p}\left\Vert F\right\Vert _{H_{p}}^{p},\text{ \ \ \ }\left( 0<p\leq
1/2\right) .  \label{2cc}
\end{equation}

However, Goginava \cite{Goginava} \ (see also \cite{tep1}) proved that the
assumption $p>1/2$ in Theorem W1 is essential. In particular, the following
is true:

\textbf{Theorem G (Goginava). }Let $0<p\leq 1/2.$ Then there exists a
martingale $F\in H_{p},$ such that%
\begin{equation*}
\sup_{n}\left\Vert \sigma _{n}F\right\Vert _{p}=+\infty .
\end{equation*}

As main result we prove that (\ref{2cc}) holds, though Theorem W1 is not
true for $0<p\leq 1/2.$

The results for summability of Fejér means of Walsh-Fourier series can be
found in \cite{FF, Fu}, \cite{GoAMH, gog8}, \cite{PS}, \cite{Sc, Si2}, \cite%
{tep2, We5}.

\section{Definitions and Notations}

Let $\mathbb{N}_{+}$ denote the set of the positive integers, $\mathbb{N}:=%
\mathbb{N}_{+}\cup \{0\}.$ Denote by $Z_{2}$ the discrete cyclic group of
order 2, that is $Z_{2}:=\{0,1\},$ where the group operation is the modulo 2
addition and every subset is open. The Haar measure on $Z_{2}$ is given so
that the measure of a singleton is 1/2.

Define the group $G$ as the complete direct product of the group $Z_{2}$
with the product of the discrete topologies of $Z_{2}$`s.

The elements of $G$ represented by sequences
\begin{equation*}
x:=(x_{0},x_{1},...,x_{j},...)\qquad \left( x_{k}=0,1\right) .
\end{equation*}

It is easy to give a base for the neighborhood of $G$
\begin{equation*}
I_{0}\left( x\right) :=G,
\end{equation*}%
\begin{equation*}
I_{n}(x):=\{y\in G\mid y_{0}=x_{0},...,y_{n-1}=x_{n-1}\}\text{ }(x\in G,%
\text{ }n\in \mathbb{N}).
\end{equation*}%
Denote $I_{n}:=I_{n}\left( 0\right) $ for $n\in \mathbb{N}$ and $\overline{%
I_{n}}:=G$ $\backslash $ $I_{n}$.

Let

\begin{equation*}
e_{n}:=\left( 0,...,0,x_{n}=1,0,...\right) \in G\qquad \left( n\in \mathbb{N}%
\right) .
\end{equation*}

It is evident
\begin{equation}
\overline{I_{M}}=\left( \overset{M-2}{\underset{k=0}{\bigcup }}\overset{M-1}{%
\underset{l=k+1}{\bigcup }}I_{l+1}\left( e_{k}+e_{l}\right) \right) \bigcup
\left( \underset{k=0}{\bigcup\limits^{M-1}}I_{M}\left( e_{k}\right) \right) .
\label{2}
\end{equation}

If $n\in \mathbb{N},$ then every $n$ can be uniquely expressed as $%
n=\sum_{k=0}^{\infty }n_{j}2^{j}$ where $n_{j}\in Z_{2}$ $~(j\in \mathbb{N})$
and only a finite number of $n_{j}`$s differ from zero. Let $\left\vert
n\right\vert :=\max $ $\{j\in \mathbb{N},$ $n_{j}\neq 0\},$ that is $%
2^{\left\vert n\right\vert }\leq n\leq 2^{\left\vert n\right\vert +1}.$

Define the variation of an $n\in \mathbb{N}$ with binary coefficients $%
\left( n_{k},\text{ }k\in \mathbb{N}\right) $ by

\begin{equation*}
V\left( n\right) =n_{0}+\overset{\infty }{\underset{k=1}{\sum }}\left|
n_{k}-n_{k-1}\right| .
\end{equation*}

For $n=\sum_{i=1}^{s}2^{n_{i}},$ $n_{1}<n_{2}<...<n_{s}$ we denote
\begin{equation*}
n^{\left( i\right) }=2^{n_{1}}+...+2^{n_{i-1}},\text{ \ }i=2,...,s
\end{equation*}%
and
\begin{equation*}
\mathbb{A}_{0,2}=\left\{ n\in \mathbb{N}:\text{ }n=2^{0}+2^{2}+%
\sum_{i=3}^{s_{n}}2^{n_{i}}\right\} .
\end{equation*}

For every $n\in \mathbb{N}$ there exists numbers $0\leq l_{1}\leq m_{1}\leq
l_{2}-2<l_{2}\leq m_{2}\leq ...\leq l_{s}-2<l_{s}\leq m_{s}$ such that it
can be uniquely expressed as $n=\sum_{i=1}^{s}\sum_{k=l_{i}}^{m_{i}}2^{k},$
where $s$ is depending only on $n.$ It is easy to show that $s\leq V\left(
n\right) \leq 2s+1.$

Define $k$-th Rademacher functions as
\begin{equation*}
r_{k}\left( x\right) :=\left( -1\right) ^{x_{k}}\text{\qquad }\left( \text{ }%
x\in G,\text{ }k\in \mathbb{N}\right) .
\end{equation*}

Now, define the Walsh system $w:=(w_{n}:n\in \mathbb{N})$ on $G$ in the
following way:
\begin{equation*}
w_{n}(x):=\underset{k=0}{\overset{\infty }{\prod }}r_{k}^{n_{k}}\left(
x\right) =r_{\left\vert n\right\vert }\left( x\right) \left( -1\right) ^{%
\underset{k=0}{\overset{\left\vert n\right\vert -1}{\sum }}n_{k}x_{k}}\text{%
\qquad }\left( n\in \mathbb{N}\right) .
\end{equation*}

The Walsh system is orthonormal and complete in $L_{2}\left( G\right) $ (see
\cite{sws}).

If $f\in L_{1}\left( G\right) $ we can establish the Fourier coefficients,
the partial sums of the Fourier series, the Fejér means, the Dirichlet and
Fejér kernels in the usual manner:
\begin{eqnarray*}
\widehat{f}\left( k\right) &:&=\int_{G}fw_{k}d\mu \,\,\,\,\qquad \ \ \left(
k\in \mathbb{N}\right) , \\
S_{n}f &:&=\sum_{k=0}^{n-1}\widehat{f}\left( k\right) w_{k}\text{ \qquad }%
\left( n\in \mathbb{N}_{+},S_{0}f:=0\right) , \\
\qquad \sigma _{n}f &:&=\frac{1}{n}\sum_{k=1}^{n}S_{k}f\text{\qquad\ \ \ }%
\left( \text{ }n\in \mathbb{N}_{+}\text{ }\right) , \\
D_{n} &:&=\sum_{k=0}^{n-1}w_{k\text{ }}\,\,\qquad \ \ \ \ \ \ \,\left( n\in
\mathbb{N}_{+}\right) , \\
K_{n} &:&=\frac{1}{n}\overset{n}{\underset{k=1}{\sum }}D_{k}\text{ \qquad
\thinspace }\left( \text{ }n\in \mathbb{N}_{+}\text{ }\right) .
\end{eqnarray*}

Recall that
\begin{equation}
D_{2^{n}}\left( x\right) =\left\{
\begin{array}{l}
2^{n},\,\text{\ \ \ }\,\text{if\thinspace \thinspace \thinspace }x\in I_{n},
\\
0,\,\ \ \ \ \ \,\text{if}\,\,x\notin I_{n}.%
\end{array}%
\right.  \label{1dn}
\end{equation}

Let $n\in \mathbb{N}.$ Then (see \cite{sws})

\begin{equation}
\frac{1}{8}V\left( n\right) \leq \left\Vert D_{n}\right\Vert _{1}\leq
V\left( n\right) .  \label{22b}
\end{equation}%
It is well-known that (see \cite{G-E-S, sws})
\begin{equation}
\sup_{n}\int_{G}\left\vert K_{n}\left( x\right) \right\vert d\mu \left(
x\right) \leq c<\infty  \label{4}
\end{equation}%
and
\begin{equation}
nK_{n}=\sum_{r=1}^{s}\left( \underset{j=r+1}{\overset{s}{\prod }}%
w_{2^{n_{j}}}\right) 2^{n_{r}}K_{2^{n_{r}}}+\sum_{t=2}^{s}\left( \underset{%
j=t+1}{\overset{s}{\prod }}w_{2^{n_{j}}}\right) n^{\left( t\right)
}D_{2^{n_{t}}},  \label{4b}
\end{equation}%
for $n=\sum_{i=1}^{s}2^{n_{i}},$ $n_{1}<n_{2}<...<n_{s}$.

For $2^{n}$-th Fejér kernel we have following equality:
\begin{equation}
K_{2^{n}}\left( x\right) =\left\{
\begin{array}{c}
\text{ }2^{t-1},\text{\ \ \ \ \ \ \ \ if }x\in I_{n}\left( e_{t}\right) , \\
\frac{2^{n}+1}{2},\text{ \ \ \ if \ }x\in I_{n},\text{\ } \\
0,\text{ \ \ \ \ \ \ \ otherwise.\ \ }%
\end{array}%
\right.  \label{5a}
\end{equation}%
for $n>t$ and $t,n\in \mathbb{N},$ (see \cite{gat}).

The norm (or quasinorm) of the space $L_{p}(G)$ is defined by \qquad

\begin{equation*}
\left\Vert f\right\Vert _{p}:=\left( \int_{G}\left\vert f\right\vert
^{p}d\mu \right) ^{1/p}\qquad \left( 0<p<\infty \right) .
\end{equation*}

The space $L_{p,\infty }\left( G\right) $ consists of all measurable
functions $f$ \ for which

\begin{equation*}
\left\Vert f\right\Vert _{L_{p,\infty }}:=\underset{\lambda >0}{\sup }%
\lambda \mu \left( f>\lambda \right) ^{1/p}<+\infty .
\end{equation*}

The $\sigma -$algebra generated by the intervals $\left\{ I_{n}\left(
x\right) :x\in G\right\} $ will be denoted by $\digamma _{n}\left( n\in
\mathbb{N}\right) .$ The conditional expectation operators relative to $%
\digamma _{n}\left( n\in \mathbb{N}\right) $ are denoted by $E_{n}.$

A sequence $F=\left( F_{n},\text{ }n\in \mathbb{N}\right) $ of functions $%
F_{n}\in L_{1}\left( G\right) $ is said to be a dyadic martingale if (for
details see e.g. \cite{sws})

$\left( i\right) $ $F_{n}$ is $\digamma _{n}$ measurable for all $n\in
\mathbb{N},$

$\left( ii\right) $ $E_{n}F_{m}=F_{n}$ for all $n\leq m.$

The maximal function of a martingale $F$ is defined by

\begin{equation*}
F^{\ast }=\sup_{n\in \mathbb{N}}\left\vert F_{n}\right\vert .
\end{equation*}

In case $f\in L_{1}\left( G\right) ,$ the maximal functions are also be
given by

\begin{equation*}
f^{\ast }\left( x\right) =\sup\limits_{n\in \mathbb{N}}\frac{1}{\mu \left(
I_{n}\left( x\right) \right) }\left\vert \int\limits_{I_{n}\left( x\right)
}f\left( u\right) d\mu \left( u\right) \right\vert .
\end{equation*}

For $0<p<\infty $ the Hardy martingale spaces $H_{p}$ $\left( G\right) $
consist of all martingale for which

\begin{equation*}
\left\Vert F\right\Vert _{H_{p}}:=\left\Vert F^{\ast }\right\Vert
_{p}<\infty .
\end{equation*}

A bounded measurable function $a$ is p-atom, if there exist a dyadic
interval $I$, such that \qquad
\begin{equation*}
\left\{
\begin{array}{l}
a)\qquad \int_{I}ad\mu =0, \\
b)\ \qquad \left\| a\right\| _{\infty }\leq \mu \left( I\right) ^{-1/p}, \\
c)\qquad \text{supp}\left( a\right) \subset I.\qquad%
\end{array}
\right.
\end{equation*}

The dyadic Hardy martingale spaces $H_{p}$ $\left( G\right) $ for $0<p\leq 1$
have an atomic characterization. Namely the following theorem is true (see
\cite{S} and \cite{We5}):

\textbf{Theorem W2}: A martingale $F=\left( F_{n},\text{ }n\in \mathbb{N}%
\right) $ is in $H_{p}\left( 0<p\leq 1\right) $ if and only if there exists
a sequence $\left( a_{k},\text{ }k\in \mathbb{N}\right) $ of p-atoms and a
sequence $\left( \mu _{k},\text{ }k\in \mathbb{N}\right) $ of a real numbers
such that for every $n\in \mathbb{N}$

\begin{equation}
\qquad \sum_{k=0}^{\infty }\mu _{k}S_{2^{n}}a_{k}=F_{n}  \label{2A}
\end{equation}%
and

\begin{equation*}
\qquad \sum_{k=0}^{\infty }\left\vert \mu _{k}\right\vert ^{p}<\infty ,
\end{equation*}%
Moreover, $\left\Vert F\right\Vert _{H_{p}}\backsim \inf \left(
\sum_{k=0}^{\infty }\left\vert \mu _{k}\right\vert ^{p}\right) ^{1/p},$
where the infimum is taken over all decomposition of $F$ of the form (\ref%
{2A}).

It is easy to check that for every martingale $F=\left( F_{n},n\in \mathbb{N}%
\right) $ and every $k\in \mathbb{N}$ the limit

\begin{equation}
\widehat{F}\left( k\right) :=\lim_{n\rightarrow \infty }\int_{G}F_{n}\left(
x\right) w_{k}\left( x\right) d\mu \left( x\right)  \label{3a}
\end{equation}%
exists and it is called the $k$-th Walsh-Fourier coefficients of $F.$

If $F:=$ $\left( E_{n}f:n\in \mathbb{N}\right) $ is regular martingale
generated by $f\in L_{1}\left( G\right) ,$ then \qquad \qquad \qquad \qquad

\begin{equation*}
\widehat{F}\left( k\right) =\int_{G}f\left( x\right) w_{k}\left( x\right)
d\mu \left( x\right) =\widehat{f}\left( k\right) ,\text{ }k\in \mathbb{N}.
\end{equation*}%
For the martingale
\begin{equation*}
F=\overset{\infty }{\underset{n=0}{\sum }}\left( F_{n}-F_{n-1}\right)
\end{equation*}%
the conjugate transforms are defined as
\begin{equation*}
\widetilde{F^{\left( t\right) }}=\overset{\infty }{\underset{n=0}{\sum }}%
r_{n}\left( t\right) \left( F_{n}-F_{n-1}\right) ,
\end{equation*}%
where $t\in G$ is fixed. Note that $\widetilde{F^{\left( 0\right) }}=F.$ As
is well known (see \cite{We3})
\begin{equation}
\left\Vert \widetilde{F^{\left( t\right) }}\right\Vert _{H_{p}}=\left\Vert
F\right\Vert _{H_{p}}  \label{5.1}
\end{equation}%
and
\begin{equation}
\left\Vert F\right\Vert _{H_{p}}^{p}\sim \int_{G}\left\Vert \widetilde{%
F^{\left( t\right) }}\right\Vert _{p}^{p}dt.  \label{5.2}
\end{equation}

\section{Formulation of Main Results}

\begin{theorem}
a) Let $0<p\leq 1/2.$ Then there exists absolute constant $c_{p},$ depending
only on $p,$ such that
\end{theorem}

\begin{equation*}
\frac{1}{\log ^{\left[ 1/2+p\right] }n}\overset{n}{\underset{m=1}{\sum }}%
\frac{\left\Vert \sigma _{m}F\right\Vert _{H_{p}}^{p}}{m^{2-2p}}\leq
c_{p}\left\Vert F\right\Vert _{H_{p}}^{p}.
\end{equation*}%
\textit{b) Let }$0<p<1/2$\textit{\ and }$\Phi :\mathbb{N}_{+}\rightarrow
\lbrack 1,$\textit{\ }$\infty )$\textit{\ is any nondecreasing function,
satisfying the conditions }$\Phi \left( n\right) \uparrow \infty $\textit{\
and }%
\begin{equation*}
\overline{\underset{k\rightarrow \infty }{\lim }}\frac{2^{k\left(
2-2p\right) }}{\Phi \left( 2^{k}\right) }=\infty .
\end{equation*}

\textit{Then there exists a martingale }$F\in H_{p},$\textit{\ such that }%
\begin{equation*}
\underset{m=1}{\overset{\infty }{\sum }}\frac{\left\Vert \sigma
_{m}F\right\Vert _{L_{p,\infty }}^{p}}{\Phi \left( m\right) }=\infty .
\end{equation*}

\begin{theorem}
Let $F\in H_{1/2}.$ Then
\begin{equation*}
\underset{n\in \mathbf{%
\mathbb{N}
}_{+}}{\sup }\underset{\left\Vert F\right\Vert _{H_{p}}\leq 1}{\sup }\frac{1%
}{n}\underset{m=1}{\overset{n}{\sum }}\left\Vert \sigma _{m}F\right\Vert
_{1/2}^{1/2}=\infty .
\end{equation*}
\end{theorem}

\section{AUXILIARY PROPOSITIONS}

\begin{lemma}
\cite{We1} \label{lemma(Weisz)} \textit{Suppose that an operator }$T$\textit{%
\ is sublinear and for some }$0<p\leq 1$
\end{lemma}

\begin{equation*}
\int\limits_{\overset{-}{I}}\left\vert Ta\right\vert ^{p}d\mu \leq
c_{p}<\infty ,
\end{equation*}%
\textit{for every }$p$\textit{-atom }$a$\textit{, where }$I$\textit{\ denote
the support of the atom. If }$T$\textit{\ is bounded from }$L_{\infty \text{
}}$\textit{\ to }$L_{\infty },$\textit{\ then }%
\begin{equation*}
\left\Vert TF\right\Vert _{p}\leq c_{p}\left\Vert F\right\Vert _{H_{p}}.
\end{equation*}

\begin{lemma}
\cite{GoSzeged} \label{lemma(Goginava)} \textit{Let }$x\in I_{l+1}\left(
e_{k}+e_{l}\right) ,$\textit{\ }$k=0,...,M-2,$\textit{\ }$l=0,...,M-1.$%
\textit{\ Then }%
\begin{equation*}
\int_{I_{M}}\left\vert K_{n}\left( x+t\right) \right\vert d\mu \left(
t\right) \leq \frac{c2^{l+k}}{n2^{M}},\text{ \ \ \ for \ }n>2^{M}.
\end{equation*}%
\
\end{lemma}

\textit{Let }$x\in I_{M}\left( e_{k}\right) ,$\textit{\ }$m=0,...,M-1.$%
\textit{\ Then }%
\begin{equation*}
\int_{I_{M}}\left\vert K_{n}\left( x+t\right) \right\vert d\mu \left(
t\right) \leq \frac{c2^{k}}{2^{M}},\text{ \ \ \textit{for}}\mathit{\ \ }\
n>2^{M}.
\end{equation*}

The next estimation of Fejér kernel is generalized of Lemma 2 in \cite{BGG}
(see also \cite{BGG2}).

\begin{lemma}
\label{lemma(Tephnadze)} Let $n=\sum_{i=1}^{s}\sum_{k=l_{i}}^{m_{i}}2^{k},$
\ where $0\leq l_{1}\leq m_{1}\leq l_{2}-2<l_{2}\leq m_{2}\leq ...\leq
l_{s}-2<l_{s}\leq m_{s}.$ Then
\begin{equation*}
n\left\vert K_{n}\left( x\right) \right\vert \geq \frac{2^{2l_{i}}}{16},%
\text{ \ \ for \ \ }x\in I_{l_{i}+1}\left( e_{l_{i}-1}+e_{l_{i}}\right) .
\end{equation*}
\end{lemma}

\textbf{Proof:} From (\ref{4b}) we have%
\begin{eqnarray*}
nK_{n} &=&\sum_{r=1}^{s}\sum_{k=l_{r}}^{m_{r}}\left( \underset{j=r+1}{%
\overset{s}{\prod }}\underset{q=l_{j}}{\overset{m_{j}}{\prod }}w_{2^{q}}%
\underset{j=k+1}{\overset{m_{r}}{\prod }}w_{2^{j}}\right) 2^{k}K_{2^{k}} \\
&&+\sum_{r=1}^{s}\sum_{k=l_{r}}^{m_{r}}\left( \underset{j=r+1}{\overset{s}{%
\prod }}\underset{q=l_{j}}{\overset{m_{j}}{\prod }}w_{2^{q}}\underset{j=k+1}{%
\overset{m_{r}}{\prod }}w_{2^{j}}\right) \left(
\sum_{t=1}^{r-1}\sum_{q=l_{t}}^{m_{t}}2^{q}+\sum_{q=l_{r}}^{k-1}2^{q}\right)
D_{2^{k}}.
\end{eqnarray*}

Let $x\in I_{l_{i}+1}\left( e_{l_{i}-1}+e_{l_{i}}\right) .$ Then

\begin{equation*}
n\left\vert K_{n}\right\vert \geq \left\vert
2^{l_{i}}K_{2^{l_{i}}}\right\vert
-\sum_{r=1}^{i-1}\sum_{k=l_{r}}^{m_{r}}\left\vert 2^{k}K_{2^{k}}\right\vert
-\sum_{r=1}^{i-1}\sum_{k=l_{r}}^{m_{r}}\left\vert 2^{k}D_{2^{k}}\right\vert
=I-II-III.
\end{equation*}

Using (\ref{5a}) we have

\begin{equation}
I=\left\vert 2^{l_{i}}K_{2^{l_{i}}}\left( x\right) \right\vert =\frac{%
2^{2l_{i}}}{4}.  \label{10.0}
\end{equation}

Since $m_{i-1}\leq l_{i}-2$ we obtain that%
\begin{equation}
II\leq \sum_{n=0}^{l_{i}-2}\left\vert 2^{n}K_{2^{n}}\left( x\right)
\right\vert \leq \sum_{n=0}^{l_{i}-2}2^{n}\frac{\left( 2^{n}+1\right) }{2}%
\leq \frac{2^{2l_{i}}}{24}+\frac{2^{l_{i}}}{4}-\frac{2}{3}.  \label{10.1}
\end{equation}

For $III$ we have
\begin{equation}
III\leq \sum_{k=0}^{l_{i}-2}\left\vert 2^{k}D_{2^{k}}\left( x\right)
\right\vert \leq \sum_{k=0}^{l_{i}-2}4^{k}=\frac{2^{2l_{i}}}{12}-\frac{1}{3}.
\label{10.2}
\end{equation}

Combining (\ref{10.0}-\ref{10.2}) we have
\begin{equation}
n\left\vert K_{n}\left( x\right) \right\vert \geq I-II-III\geq \frac{%
2^{2l_{i}}}{8}-\frac{2^{l_{i}}}{4}+1.  \label{10.3}
\end{equation}

Suppose that $l_{i}\geq 2$. Then
\begin{equation*}
n\left\vert K_{n}\left( x\right) \right\vert \geq \frac{2^{2l_{i}}}{8}-\frac{%
2^{2l_{i}}}{16}\geq \frac{2^{2l_{i}}}{16}.
\end{equation*}

If $l_{i}=0$ or $l_{i}=1$ using (\ref{10.3}) we have
\begin{equation*}
n\left\vert K_{n}\left( x\right) \right\vert \geq \frac{7}{8}\geq \frac{%
2^{2l_{i}}}{16}.
\end{equation*}

Lemma 3 is proved.

\section{Proof of the Theorems}

\textbf{Proof of Theorem 1. }Suppose that
\begin{equation*}
\frac{1}{\log ^{\left[ 1/2+p\right] }n}\overset{n}{\underset{m=1}{\sum }}%
\frac{\left\Vert \sigma _{m}F\right\Vert _{p}^{p}}{m^{2-2p}}\leq
c_{p}\left\Vert F\right\Vert _{H_{p}}^{p}.
\end{equation*}

Combining (\ref{5.1}) and (\ref{5.2}) we have
\begin{eqnarray}
&&\frac{1}{\log ^{\left[ 1/2+p\right] }n}\overset{n}{\underset{m=1}{\sum }}%
\frac{\left\Vert \sigma _{m}F\right\Vert _{H_{p}}^{p}}{m^{2-2p}}  \label{5.3}
\\
&=&\frac{1}{\log ^{\left[ 1/2+p\right] }n}\overset{n}{\underset{m=1}{\sum }}%
\frac{\int_{G}\left\Vert \widetilde{\sigma _{m}F^{\left( t\right) }}%
\right\Vert _{p}^{p}dt}{m^{2-2p}}  \notag \\
&=&\frac{1}{\log ^{\left[ 1/2+p\right] }n}\overset{n}{\underset{m=1}{\sum }}%
\frac{\int_{G}\left\Vert \sigma _{m}\widetilde{F^{\left( t\right) }}%
\right\Vert _{p}^{p}dt}{m^{2-2p}}  \notag \\
&=&\int_{G}\frac{1}{\log ^{\left[ 1/2+p\right] }n}\overset{n}{\underset{m=1}{%
\sum }}\frac{\left\Vert \sigma _{m}\widetilde{F^{\left( t\right) }}%
\right\Vert _{p}^{p}}{m^{2-2p}}dt  \notag \\
&\leq &c_{p}\int_{G}\left\Vert \widetilde{F^{\left( t\right) }}\right\Vert
_{H_{p}}^{p}dt\sim c_{p}\int_{G}\left\Vert F\right\Vert
_{H_{p}}^{p}dt=c_{p}\left\Vert F\right\Vert _{H_{p}}^{p}.  \notag
\end{eqnarray}

By Lemma \ref{lemma(Weisz)} and (\ref{5.3}) the proof of theorem 1 will be
complete, if we show that

\begin{equation*}
\frac{1}{\log ^{\left[ 1/2+p\right] }n}\overset{n}{\underset{m=1}{\sum }}%
\frac{\left\Vert \sigma _{m}a\right\Vert _{p}^{p}}{m^{2-2p}}\leq c<\infty ,%
\text{ \ \ \ }m=1,2,...
\end{equation*}%
for every p-atom $a.$ We may assume that $a$ be an arbitrary p-atom with
support$\ I$, $\mu \left( I\right) =2^{-M}$ and $I=I_{M}.$ It is easy to see
that $\sigma _{n}\left( a\right) =0,$ when $n\leq 2^{M}.$ Therefore we can
suppose that $n>2^{M}.$

Let $x\in I_{M}.$ Since $\sigma _{n}$ is bounded from $L_{\infty }$ to $%
L_{\infty }$ (the boundedness follows from (\ref{4})) and $\left\Vert
a\right\Vert _{\infty }\leq 2^{M/p}$ we obtain
\begin{equation*}
\int_{I_{M}}\left\vert \sigma _{m}a\left( x\right) \right\vert ^{p}d\mu
\left( x\right) \leq \left\Vert \sigma _{m}a\right\Vert _{\infty
}^{p}/2^{M}\leq c<\infty ,\text{ }0<p\leq 1/2.
\end{equation*}%
Let $0<p\leq 1/2.$ Hence
\begin{equation*}
\frac{1}{\log ^{\left[ 1/2+p\right] }n}\overset{n}{\underset{m=1}{\sum }}%
\frac{\int_{I_{M}}\left\vert \sigma _{m}a\left( x\right) \right\vert
^{p}d\mu \left( x\right) }{m^{2-2p}}\leq \frac{c}{\log ^{\left[ 1/2+p\right]
}n}\overset{n}{\underset{m=1}{\sum }}\frac{1}{m^{2-2p}}\leq c<\infty .
\end{equation*}

It is easy to show that
\begin{equation*}
\left\vert \sigma _{m}a\left( x\right) \right\vert \leq
\int_{I_{M}}\left\vert a\left( t\right) \right\vert \left\vert K_{m}\left(
x+t\right) \right\vert d\mu \left( t\right) \leq
2^{M/p}\int_{I_{M}}\left\vert K_{m}\left( x+t\right) \right\vert d\mu \left(
t\right) .
\end{equation*}

From Lemma \ref{lemma(Goginava)} we get
\begin{equation}
\left\vert \sigma _{m}a\left( x\right) \right\vert \leq \frac{%
c2^{k+l}2^{M\left( 1/p-1\right) }}{m},\text{ \ \ }x\in I_{l+1}\left(
e_{k}+e_{l}\right) ,\,0\leq k<l<M,  \label{12}
\end{equation}%
and

\begin{equation}
\left\vert \sigma _{m}a\left( x\right) \right\vert \leq c2^{M\left(
1/p-1\right) }2^{k},\text{ \ }x\in I_{M}\left( e_{k}\right) ,\,0\leq k<M.
\label{12a}
\end{equation}

Combining (\ref{2}) and (\ref{12}-\ref{12a}) we obtain
\begin{eqnarray}
&&\int_{\overline{I_{M}}}\left\vert \sigma _{m}a\left( x\right) \right\vert
^{p}d\mu \left( x\right)  \label{7.2} \\
&=&\overset{M-2}{\underset{k=0}{\sum }}\overset{M-1}{\underset{l=k+1}{\sum }}%
\int_{I_{l+1}\left( e_{k}+e_{l}\right) }\left\vert \sigma _{m}a\left(
x\right) \right\vert ^{p}d\mu \left( x\right) +\overset{M-1}{\underset{k=0}{%
\sum }}\int_{I_{M}\left( e_{k}\right) }\left\vert \sigma _{m}a\left(
x\right) \right\vert ^{p}d\mu \left( x\right)  \notag \\
&\leq &c\overset{M-2}{\underset{k=0}{\sum }}\overset{M-1}{\underset{l=k+1}{%
\sum }}\frac{1}{2^{l}}\frac{2^{p\left( k+l\right) }2^{M\left( 1-p\right) }}{%
m^{p}}+c\overset{M-1}{\underset{k=0}{\sum }}\frac{1}{2^{M}}2^{M\left(
1-p\right) }2^{pk}  \notag \\
&\leq &\frac{c2^{M\left( 1-p\right) }}{m^{p}}\overset{M-2}{\underset{k=0}{%
\sum }}\overset{M-1}{\underset{l=k+1}{\sum }}\frac{2^{p\left( k+l\right) }}{%
2^{l}}+c\overset{M-1}{\underset{k=0}{\sum }}\frac{2^{pk}}{2^{pM}}  \notag \\
&\leq &\frac{c2^{M\left( 1-p\right) }M^{\left[ 1/2+p\right] }}{m^{p}}+c.
\notag
\end{eqnarray}

\bigskip Hence

\begin{eqnarray*}
&&\frac{1}{\log ^{\left[ 1/2+p\right] }n}\overset{n}{\underset{m=2^{M}+1}{%
\sum }}\frac{\int_{\overline{I_{M}}}\left\vert \sigma _{m}a\left( x\right)
\right\vert ^{p}d\mu \left( x\right) }{m^{2-2p}} \\
&\leq &\frac{1}{\log ^{\left[ 1/2+p\right] }n}\left( \overset{n}{\underset{%
m=2^{M}+1}{\sum }}\frac{c2^{M\left( 1-p\right) }M^{\left[ 1/2+p\right] }}{%
m^{2-p}}+\overset{n}{\underset{m=2^{M}+1}{\sum }}\frac{c}{m^{2-2p}}\right)
<c<\infty .
\end{eqnarray*}

Which complete the proof of first part of Theorem 1.

\bigskip Now we proof second part of Theorem 1. Let $\Phi \left( n\right) $
is nondecreasing function, satisfying condition
\begin{equation}
\underset{k\rightarrow \infty }{\lim }\frac{2^{\left( \left\vert
n_{k}\right\vert +1\right) \left( 2-2p\right) }}{\Phi \left( 2^{\left\vert
n_{k}\right\vert +1}\right) }=\infty .  \label{12j}
\end{equation}

Under condition (\ref{12j}), there exists a sequence $\left\{ \alpha _{k}:%
\text{ }k\geq 0\right\} \subset \left\{ n_{k}:\text{ }k\geq 0\right\} $ such
that
\begin{equation}
\left\vert \alpha _{k}\right\vert \geq 2,\text{ \ \ \ for all \ }k\geq 0
\label{122}
\end{equation}%
and
\begin{equation}
\sum_{\eta =0}^{\infty }\frac{\Phi ^{1/2}\left( 2^{\left\vert \alpha _{\eta
}\right\vert +1}\right) }{2^{\left\vert \alpha _{\eta }\right\vert \left(
1-p\right) }}=2^{1-p}\sum_{\eta =0}^{\infty }\frac{\Phi ^{1/2}\left(
2^{\left\vert \alpha _{\eta }\right\vert +1}\right) }{2^{\left( \left\vert
\alpha _{\eta }\right\vert +1\right) \left( 1-p\right) }}<c<\infty .
\label{121}
\end{equation}

Let \qquad
\begin{equation*}
F_{n}=\sum_{\left\{ k:\text{ }\left\vert \alpha _{k}\right\vert <n\right\}
}\lambda _{k}a_{k},\text{ \ }
\end{equation*}%
where
\begin{equation*}
\lambda _{k}=\frac{\Phi ^{1/2p}\left( 2^{\left\vert \alpha _{k}\right\vert
+1}\right) }{2^{\left( \left\vert \alpha _{k}\right\vert \right) \left(
1/p-1\right) }}
\end{equation*}%
and

\begin{equation*}
a_{k}=2^{\left\vert \alpha _{k}\right\vert \left( 1/p-1\right) }\left(
D_{2^{\left\vert \alpha _{k}\right\vert +1}}-D_{2^{\left\vert \alpha
_{k}\right\vert }}\right) .
\end{equation*}

It is easy to show that the martingale $\,F=\left( F_{n},\text{ }n\in
\mathbb{N}\right) \in H_{p}.$

Indeed, since

\begin{equation*}
S_{2^{n}}a_{k}=\left\{
\begin{array}{l}
a_{k},\text{ \quad }\left\vert \alpha _{k}\right\vert <n, \\
0,\text{ \ \ \ \ \ }\left\vert \alpha _{k}\right\vert \geq n,%
\end{array}%
\right.
\end{equation*}

\begin{eqnarray*}
\text{supp}(a_{k}) &=&I_{\left\vert \alpha _{k}\right\vert }, \\
\int_{I_{\left\vert \alpha _{k}\right\vert }}a_{k}d\mu &=&0
\end{eqnarray*}%
and

\begin{equation*}
\left\Vert a_{k}\right\Vert _{\infty }\leq 2^{\left\vert \alpha
_{k}\right\vert /p}=\left( \text{supp }a_{k}\right) ^{-1/p}
\end{equation*}%
if we apply Theorem W2 and (\ref{121}) we conclude that $F\in H_{p}.$

It is easy to show that

\begin{eqnarray}
&&\widehat{F}(j)  \label{6aa} \\
&=&\left\{
\begin{array}{l}
\Phi ^{1/2p}\left( 2^{\left\vert \alpha _{k}\right\vert +1}\right) ,\,\ \,%
\text{ if \thinspace \thinspace }j\in \left\{ 2^{\left\vert \alpha
_{k}\right\vert },...,2^{\left\vert \alpha _{k}\right\vert +1}-1\right\} ,%
\text{ }k=0,1,2..., \\
0\text{ },\text{ \thinspace \qquad \thinspace\ \ \ \ \ \ \ \ \ \ \ \
\thinspace if \thinspace \thinspace \thinspace }j\notin
\bigcup\limits_{k=0}^{\infty }\left\{ 2^{\left\vert \alpha _{k}\right\vert
},...,2^{\left\vert \alpha _{k}\right\vert +1}-1\right\} .\text{ }%
\end{array}%
\right.  \notag
\end{eqnarray}

Let $2^{\left\vert \alpha _{k}\right\vert }<n<2^{\left\vert \alpha
_{k}\right\vert +1}.$ Using (\ref{6aa}) we can write

\begin{equation}
\sigma _{_{n}}F=\frac{1}{n}\sum_{j=1}^{2^{\left\vert \alpha _{k}\right\vert
}}S_{j}F+\frac{1}{n}\sum_{j=2^{\left\vert \alpha _{k}\right\vert
}+1}^{n}S_{j}F=III+IV.  \label{7aa}
\end{equation}

It is simple to show that%
\begin{equation}
S_{j}F=\left\{
\begin{array}{l}
0,\,\ \text{ if \thinspace \thinspace }0\leq j\leq 2^{\left\vert \alpha
_{0}\right\vert } \\
\Phi ^{1/2p}\left( 2^{\left\vert \alpha _{0}\right\vert +1}\right) \left(
D_{_{j}}-D_{2^{\left\vert \alpha _{0}\right\vert }}\right) ,\text{
\thinspace\ \thinspace if \thinspace \thinspace \thinspace }2^{\left\vert
\alpha _{0}\right\vert }<j\leq 2^{\left\vert \alpha _{0}\right\vert +1}.%
\text{ }%
\end{array}%
\right.  \label{7aaa}
\end{equation}

Suppose that $2^{\left\vert \alpha _{s}\right\vert }<j\leq 2^{\left\vert
\alpha _{s}\right\vert +1},$ for some $s=1,2,...,k.$ Then applying (\ref{6aa}%
) we have%
\begin{eqnarray}
S_{j}F &=&\sum_{v=0}^{2^{\left\vert \alpha _{s-1}\right\vert +1}-1}\widehat{F%
}(v)w_{v}+\sum_{v=2^{^{\left\vert \alpha _{s}\right\vert }}}^{j-1}\widehat{F}%
(v)w_{v}  \label{8aa} \\
&=&\sum_{\eta =0}^{s-1}\sum_{v=2^{\left\vert \alpha _{\eta }\right\vert
}}^{2^{\left\vert \alpha _{\eta }\right\vert +1}-1}\widehat{F}%
(v)w_{v}+\sum_{v=2^{\left\vert \alpha _{s}\right\vert }}^{j-1}\widehat{F}%
(v)w_{v}  \notag \\
&=&\sum_{\eta =0}^{s-1}\sum_{v=2^{\left\vert \alpha _{\eta }\right\vert
}}^{2^{\left\vert \alpha _{\eta }\right\vert +1}-1}\Phi ^{1/2p}\left(
2^{\left\vert \alpha _{\eta }\right\vert +1}\right) w_{v}+\Phi ^{1/2p}\left(
2^{\left\vert \alpha _{s}\right\vert +1}\right) \sum_{v=2^{\left\vert \alpha
_{s}\right\vert }}^{j-1}w_{v}  \notag \\
&=&\sum_{\eta =0}^{s-1}\Phi ^{1/2p}\left( 2^{\left\vert \alpha _{\eta
}\right\vert +1}\right) \left( D_{2^{\left\vert \alpha _{\eta }\right\vert
+1}}-D_{2^{\left\vert \alpha _{\eta }\right\vert }}\right)  \notag \\
&&+\Phi ^{1/2p}\left( 2^{\left\vert \alpha _{s}\right\vert +1}\right) \left(
D_{_{j}}-D_{2^{\left\vert \alpha _{s}\right\vert }}\right) .  \notag
\end{eqnarray}

Let $2^{\left\vert \alpha _{s}\right\vert +1}\leq j\leq 2^{\left\vert \alpha
_{s+1}\right\vert },$ $s=0,1,...k-1.$ Analogously of (\ref{8aa}) we get%
\begin{equation}
S_{j}F=\sum_{v=0}^{2^{\left\vert \alpha _{s}\right\vert +1}}\widehat{F}%
(v)w_{v}=\sum_{\eta =0}^{s}\Phi ^{1/2p}\left( 2^{\left\vert \alpha _{\eta
}\right\vert +1}\right) \left( D_{2^{\left\vert \alpha _{\eta }\right\vert
+1}}-D_{2^{\left\vert \alpha _{\eta }\right\vert }}\right) .  \label{10aa}
\end{equation}

Let $x\in I_{2}\left( e_{0}+e_{1}\right) .$ Since (see (\ref{1dn}) and (\ref%
{5a}))
\begin{equation}
D_{2^{n}}\left( x\right) =K_{2^{n}}\left( x\right) =0,\text{ for }n\geq 2,
\label{40}
\end{equation}%
from (\ref{122}) and (\ref{8aa}-\ref{10aa}) we obtain%
\begin{equation}
III=\frac{1}{n}\sum_{\eta =0}^{k-1}\Phi ^{1/2p}\left( 2^{\left\vert \alpha
_{\eta }\right\vert +1}\right) \sum_{v=2^{\left\vert \alpha _{\eta
}\right\vert }+1}^{2^{\left\vert \alpha _{\eta }\right\vert +1}}D_{v}\left(
x\right)  \label{9aaa}
\end{equation}%
\begin{equation*}
=\frac{1}{n}\sum_{\eta =0}^{k-1}\Phi ^{1/2p}\left( 2^{\left\vert \alpha
_{\eta }\right\vert +1}\right) \left( 2^{\left\vert \alpha _{\eta
}\right\vert +1}K_{2^{\left\vert \alpha _{\eta }\right\vert +1}}\left(
x\right) -2^{\left\vert \alpha _{\eta }\right\vert }K_{2^{\left\vert \alpha
_{\eta }\right\vert }}\left( x\right) \right) =0.
\end{equation*}

Applying (\ref{8aa}), when $s=k$ in $IV$ we have

\begin{eqnarray}
IV &=&\frac{n-2^{\left\vert \alpha _{k}\right\vert }}{n}\sum_{\eta
=0}^{k-1}\Phi ^{1/2p}\left( 2^{\left\vert \alpha _{\eta }\right\vert
+1}\right) \left( D_{2^{\left\vert \alpha _{\eta }\right\vert
+1}}-D_{2^{\left\vert \alpha _{\eta }\right\vert }}\right)  \label{9aa} \\
&&+\frac{\Phi ^{1/2p}\left( 2^{\left\vert \alpha _{k}\right\vert +1}\right)
}{n}\sum_{j=2^{_{\left\vert \alpha _{k}\right\vert }}+1}^{n}\left(
D_{_{j}}-D_{2^{\left\vert \alpha _{k}\right\vert }}\right) =IV_{1}+IV_{2}.
\notag
\end{eqnarray}

Combining (\ref{122}) and (\ref{40}) we get

\begin{equation}
IV_{1}=0,\text{ \ for \ }x\in I_{2}\left( e_{0}+e_{1}\right) .  \label{8aaaa}
\end{equation}

Let $\alpha _{k}\in \mathbb{A}_{0,2},$ $2^{\left\vert \alpha _{k}\right\vert
}<n<2^{\left\vert \alpha _{k}\right\vert +1}$ and $x\in I_{2}\left(
e_{0}+e_{1}\right) $. Since $n-2^{_{\left\vert \alpha _{k}\right\vert }}\in
\mathbb{A}_{0,2}$ and%
\begin{equation}
D_{j+2^{m}}=D_{2^{m}}+w_{_{2^{m}}}D_{j},\text{ when \thinspace \thinspace }%
j<2^{m}  \label{30}
\end{equation}%
from (\ref{4b}) and (\ref{40}) we obtain
\begin{equation}
\left\vert IV_{2}\right\vert =\frac{\Phi ^{1/2p}\left( 2^{\left\vert \alpha
_{k}\right\vert +1}\right) }{n}\left\vert \sum_{j=1}^{n-2^{_{\left\vert
\alpha _{k}\right\vert }}}\left( D_{j+2^{\left\vert \alpha _{k}\right\vert
}}\left( x\right) -D_{2^{\left\vert \alpha _{k}\right\vert }}\left( x\right)
\right) \right\vert   \label{31}
\end{equation}%
\begin{eqnarray*}
&=&\frac{\Phi ^{1/2p}\left( 2^{\left\vert \alpha _{k}\right\vert +1}\right)
}{n}\left\vert \sum_{j=1}^{n-2^{\left\vert _{\alpha _{k}}\right\vert
}}D_{_{j}}\left( x\right) \right\vert  \\
&=&\frac{\Phi ^{1/2p}\left( 2^{\left\vert \alpha _{k}\right\vert +1}\right)
}{n}\left\vert \left( n-2^{_{\left\vert \alpha _{k}\right\vert }}\right)
K_{n-2^{_{\left\vert \alpha _{k}\right\vert }}}\left( x\right) \right\vert
\geq \frac{\Phi ^{1/2p}\left( 2^{\left\vert \alpha _{k}\right\vert
+1}\right) }{2^{\left\vert \alpha _{k}\right\vert +1}}.
\end{eqnarray*}

Let $0<p<1/2$ and $n\in \mathbb{A}_{0,2}.$ Combining (\ref{7aa}-\ref{31}) we
have

\begin{eqnarray}
&&\left\Vert \sigma _{n}F\right\Vert _{L_{p,\infty }}^{p}  \label{10aaa} \\
&\geq &\frac{c_{p}\Phi ^{1/2}\left( 2^{\left\vert \alpha _{k}\right\vert
+1}\right) }{2^{p\left( \left\vert \alpha _{k}\right\vert +1\right) }}\mu
\left\{ x\in I_{2}\left( e_{0}+e_{1}\right) :\text{ }\left\vert \sigma
_{n}F\right\vert \geq \frac{c_{p}\Phi ^{1/2p}\left( 2^{\left\vert \alpha
_{k}\right\vert +1}\right) }{2^{\left\vert \alpha _{k}\right\vert +1}}%
\right\}   \notag \\
&\geq &\frac{c_{p}\Phi ^{1/2}\left( 2^{\left\vert \alpha _{k}\right\vert
+1}\right) }{2^{p\left( \left\vert \alpha _{k}\right\vert +1\right) }}\mu
\left\{ I_{2}\left( e_{0}+e_{1}\right) \right\} \geq \frac{c_{p}\Phi
^{1/2}\left( 2^{\left\vert \alpha _{k}\right\vert +1}\right) }{2^{p\left(
\left\vert \alpha _{k}\right\vert +1\right) }}.  \notag
\end{eqnarray}

Hence
\begin{eqnarray*}
&&\underset{n=1}{\overset{\infty }{\sum }}\frac{\left\Vert \sigma
_{n}F\right\Vert _{L_{p,\infty }}^{p}}{\Phi \left( n\right) }\geq \underset{%
\left\{ n\in \mathbb{A}_{0,2}:\text{ }2^{\left\vert \alpha _{k}\right\vert
}<n<2^{\left\vert \alpha _{k}\right\vert +1}\right\} }{\sum }\frac{%
\left\Vert \sigma _{n}F\right\Vert _{L_{p,\infty }}^{p}}{\Phi \left(
n\right) } \\
&\geq &\frac{1}{\Phi ^{1/2}\left( 2^{\left\vert \alpha _{k}\right\vert
+1}\right) }\underset{\left\{ n\in \mathbb{A}_{0,2}:\text{ }2^{\left\vert
\alpha _{k}\right\vert }<n<2^{\left\vert \alpha _{k}\right\vert +1}\right\} }%
{\sum }\frac{1}{2^{p\left( \left\vert \alpha _{k}\right\vert +1\right) }} \\
&\geq &\frac{c_{p}2^{\left( 1-p\right) \left( \left\vert \alpha
_{k}\right\vert +1\right) }}{\Phi ^{1/2}\left( 2^{\left\vert \alpha
_{k}\right\vert +1}\right) }\rightarrow \infty ,\text{ \qquad when \ \ }%
k\rightarrow \infty .
\end{eqnarray*}

Which complete the proof of Theorem 1.

\textbf{Proof of Theorem 2. }Let

\begin{equation*}
F_{m}=2^{m}\left( D_{2^{m}+1}-D_{2^{m}}\right) .\text{ }
\end{equation*}

It is evident
\begin{equation}
\widehat{F}_{m}\left( i\right) =\left\{
\begin{array}{l}
\text{ }2^{m},\text{ if }i=2^{m},...,2^{m+1}-1, \\
\text{ }0,\text{otherwise.}%
\end{array}%
\right.  \label{14a}
\end{equation}%
We can write
\begin{equation}
S_{i}F_{m}=\left\{
\begin{array}{l}
2^{m}\left( D_{i}-D_{2^{m}}\right) ,\text{ \ if }i=2^{m}+1,...,2^{m+1}-1, \\
\text{ }F_{m},\text{ \ if }i\geq 2^{m+1}, \\
0,\text{ \ otherwise.}%
\end{array}%
\right.  \label{14}
\end{equation}

From (\ref{1dn}) we get
\begin{eqnarray}
\left\Vert F_{m}\right\Vert _{H_{1/2}} &=&\left\Vert \sup\limits_{k\in
\mathbb{N}}\left\vert S_{2^{k}}F_{m}\right\vert \right\Vert
_{1/2}=2^{m}\left\Vert D_{2^{m}+1}-D_{2^{m}}\right\Vert _{1/2}  \label{14c}
\\
&=&2^{m}\left\Vert D_{2^{m}}\right\Vert _{1/2}\leq 1.  \notag
\end{eqnarray}

Let $0<n<2^{m}.$ Using (\ref{30}) we obtain

\begin{eqnarray}
\left\vert \sigma _{n+2^{m}}F_{m}\right\vert &=&\frac{1}{n+2^{m}}\left\vert
\overset{n+2^{m}}{\underset{j=2^{m}+1}{\sum }}S_{j}F_{m}\right\vert
\label{16b} \\
&=&\frac{1}{n+2^{m}}\left\vert 2^{m}\overset{n+2^{m}}{\underset{j=2^{m}+1}{%
\sum }}\left( D_{j}-D_{2^{m}}\right) \right\vert  \notag \\
&=&\frac{1}{n+2^{m}}\left\vert 2^{m}\overset{n}{\underset{j=1}{\sum }}\left(
D_{j+2^{m}}-D_{2^{m}}\right) \right\vert  \notag \\
&=&\frac{1}{n+2^{m}}\left\vert 2^{m}\overset{n}{\underset{j=1}{\sum }}%
D_{j}\right\vert =\frac{2^{m}}{n+2^{m}}n\left\vert K_{n}\right\vert .  \notag
\end{eqnarray}

Let $n=\sum_{i=1}^{s}\sum_{k=l_{i}}^{m_{i}}2^{k},$ \ where $0\leq l_{1}\leq
m_{1}\leq l_{2}-2<l_{2}\leq m_{2}\leq ...\leq l_{s}-2<l_{s}\leq m_{s}.$
Combining Lemma \ref{lemma(Tephnadze)} and (\ref{16b}) we write
\begin{equation*}
\left\vert \sigma _{n+2^{m}}F_{m}\left( x\right) \right\vert \geq
c2^{2l_{i}},\text{ \ \ \ for \ }x\in I_{l_{i}+1}\left(
e_{l_{i}-1}+e_{l_{i}}\right) .
\end{equation*}

It follows that
\begin{eqnarray*}
&&\int_{G}\left\vert \sigma _{n+2^{m}}F_{m}(x)\right\vert ^{1/2}d\mu \left(
x\right) \\
&\geq &\text{ }\underset{i=0}{\overset{s}{\sum }}\int_{I_{l_{i}+1}\left(
e_{l_{i}-1}+e_{l_{i}}\right) }\left\vert \sigma
_{n+2^{m}}F_{m}(x)\right\vert ^{1/2}d\mu \left( x\right) \\
&\geq &c\overset{s}{\underset{i=0}{\sum }}\frac{1}{2^{l_{i}}}2^{l_{i}}\geq
cs\geq cV\left( n\right) .
\end{eqnarray*}%
On the other hand, Fine \cite{FF} proved that
\begin{equation*}
\frac{1}{n\log n}\underset{k=1}{\overset{n}{\sum }}V\left( k\right) =\frac{1%
}{4\log 2}+o\left( 1\right) .
\end{equation*}

Hence
\begin{eqnarray*}
&&\underset{n\in \mathbf{%
\mathbb{N}
}_{+}}{\sup }\underset{\left\Vert f\right\Vert _{H_{p}}\leq 1}{\sup }\frac{1%
}{n}\underset{k=1}{\overset{n}{\sum }}\left\Vert \sigma _{k}F\right\Vert
_{1/2}^{1/2} \\
&\geq &\frac{1}{2^{m+1}}\overset{2^{m+1}-1}{\underset{k=2^{m}+1}{\sum }}%
\left\Vert \sigma _{k}F_{m}\right\Vert _{1/2}^{1/2} \\
&\geq &\frac{c}{2^{m+1}}\underset{k=2^{m}+1}{\overset{2^{m+1}-1}{\sum }}%
V\left( k-2^{m}\right) \geq \frac{c}{2^{m+1}}\underset{k=1}{\overset{2^{m}-1}%
{\sum }}V\left( k\right) \\
&\geq &c\log m\rightarrow \infty ,\text{ when \ }m\rightarrow \infty .
\end{eqnarray*}

Which complete the proof of Theorem 2.

\textbf{Acknowledgment: }The author would like to thank the referee for
helpful suggestions.

\end{document}